\documentclass[preprint,prl,aps]{revtex4}

\begin{document}

%-----------------------------------------------------------------
%\usepackage{amssymb,amsmath}
%\usepackage{graphicx}
%\usepackage{subfigure,psfig}
\newcommand{\be}{\begin{equation}}
\newcommand{\ee}{\end{equation}}
\newcommand{\ben}{\begin{eqnarray}}
\newcommand{\een}{\end{eqnarray}}
\newcommand{\n}{\nonumber  }
\newcommand{\nn}{\nonumber \\ }
\newcommand{\nd}{\noindent}
\newcommand{\p}{\partial}
%-----------------------------------------------------------------

% Full title of the paper (Capitalized)
\title{Inferring an optimal Fisher measure}

% Authors (add full first names)
\author{S.P. Flego $^2$}
\author{A. Plastino $^{1,\,3,\,5}$}
\author{A. R. Plastino$^{3,\,4}$}

% Affiliations / Addresses, add [1] after \address if there is only one affiliation
\affiliation{ $^{1}$ Universidad Nacional de La Plata, Instituto
de F\'{\i}sica (IFLP-CCT-CONICET),
C.C. 727, 1900 La Plata, Argentina \\
$^{2}$ Universidad Nacional de La Plata, Fac. de Ingenier\'{\i}a,
 1900 La Plata, Argentina \\  $^{3}$  CREG-University of La Plata-CONICET\\C.C. 727, 1900 La Plata, Argentina
 \\
$^{4}$ Instituto Carlos I de Fisica Teorica y Computacional and
Departamento de Fisica Atomica, Molecular y Nuclear, Universidad
de Granada, Granada, Spain \\ $^{5}$ Universitat de les Illes
Balears and IFISC-CSIC, 07122 Palma de Mallorca, Spain}
%

%%%%%%%%%%%%%%%%%%%%%%%%%%%%%%%%%%%\corres{plastino@fisica.unlp.edu.ar}
\begin{abstract}
\nd \nd It is well known that a suggestive  relation exists that
links Schr\"odinger's equation (SE) to the information-optimizing
principle based on Fisher's  information measure (FIM). We explore
here an approach that will allow one to infer the optimal FIM
compatible with a given amount of prior information without
explicitly solving first the associated SE. This technique is
based on the virial theorem and it provides  analytic solutions
   for the physically relevant FIM, that which is minimal subject to
the constraints posed by the prior information.

\vspace{2.cm}


%%%%%%%%%\pacs{PACS. numbers: 05.45+b, 05.30-d}
\end{abstract}

%\newpage

\maketitle

%%%%%%%%%%%%%%%%%%%%%%%%%%%%%%%%%%%%%%%%%%%%%%%%%%%%%%%%%%%%
\section{I. Introduction}

\nd Although Fisher's information measure (FIM) $I$ dates from the
20's, it has attracted  intense attention from physicists only
since the mid-90's \cite{emergent,garba,frieden}. Interest in
FIM's physical applications
\cite{frieden,frieden2,frieden3,pla7,pla11,flego,reginatto} has
been growing exponentially since the appearance of Frieden and
Soffer's seminal paper \cite{frieden3}. A very small, and
certainly not exhaustive sample is that of Refs. \cite{
hall,KSC10,FS09,U09,LADY08,2011,SAA07,N07,N06,alp2,alp3,alp4,olivares,pla5,flavia}.

\nd %%In this work we are not concerned at all with statistics but
%%%only with physical-chemical FIM applications.
The mathematical problem of extremizing $I$ under given
constraints has proved to be relevant in connection with several
scenarios related to quantum mechanics and statistical physics, as
the above cited references indicate. We intend to exhibit novel
links between $I$ and the  Schr\"odinger equation that, via the
virial theorem,  lead in natural fashion to a differential
equation for $I$. Such equation, that has an analytical solution,
encodes the available prior knowledge concerning the system at
hand in terms of adequately selected expectation values. Our
solution adds to the rather large Fisher literature a general,
explicit expression
 for that particular FIM $I_{Min}$ that arises out of any constrained
 $I-$extremization problem.
 \vskip 3mm
 \nd To better understand why this is of importance one should
 recall that Fisher's information and Shannon's entropy play
 complementary roles \cite{frieden2}. The former is convex, the
 later concave. When one grows, the other diminishes, etc. The
 associated Shannon's MaxEnt problem has as its solution, always, an
 exponential form that contains those physical quantities whose
 mean values are a priori known. The solution of the FIM
 minimization problem is instead a Schr\"odinger like differential
 equation\cite{pla7,pla11,reginatto}, whose solutions exhibit a panoply
 of different mathematical forms. We provide here special definite forms for
 $I_{Min}$, expressed in terms of those mean values that are a priori
 known, filling thus a gap in the literature of the physics of
 information. We review basic materials below and start with our
 presentation in Section III-IV.  We give illustrative examples in Section V and
 discuss particular issues in Section VI. Our conclusions are given in Section VII.

%%%%%%%%%%%%%%%%%%%%%%%%%%%%%%%%%%%%%%%%%%%%%%%%%%%%%%%

\section{II. Review of basic ideas}

\nd We briefly review here the formalism developed in Ref.
\cite{pla7} (see also \cite{pla11}). Consider a system that is
specified by a physical parameter $\theta$ and let $f(x,\theta)$
describe the normalized probability distribution function (PDF)
for this parameter. If an observer were to make a measurement of
$x$ and had to best infer $\theta$ from such  measurement, calling
the resulting estimate $\tilde \theta=\tilde \theta(x)$, one might
well wonder how well $\theta$ could be determined. Estimation
theory~\cite{frieden2} asserts that the {\it best possible
estimator} $\tilde \theta( x)$, after a very large number of
$x$-samples is examined, suffers a mean-square error $e^2$ from
$\theta$ obeying the rule $Ie^2=1$, where the Fisher information
measure (FIM) $I$, a functional of the PDF, reads
 \be \label{eq.1-1} I \,=\,\int ~dx ~f(x,\theta)
\left\{\frac{\partial ~ }{\partial
\theta}~\ln{[f(x,\theta)]}\right\}^2.\ee
Any other estimator must have a larger mean-square error (all estimators must be unbiased, i.e., satisfy
$ \langle \tilde \theta({\bf x}) \rangle=\,\theta \label{unbias}$).
Thus, FIM has a lower bound. No matter what the parameter $\xi$ of the system might
be, $I$ has to obey
\be \label{rao} I\,e_\xi^2\,\ge \,1,\ee the celebrated Cramer--Rao bound \cite{frieden2}. The particular instance of translational families merits a word.
They are mono-parametric distribution families of the form
$f(x,\theta)=f(x-\theta),$  known up to the shift parameter $\theta$. All
family members exhibit identical shape. After introducing the amplitudes $\psi$ such that $f(x)=\psi(x)^2$, FIM adopts the simpler aspect \cite{frieden3}
\be \label{eq.1-2} I \,=\,\int ~dx ~f(x)\left\{\frac{\partial ~
}{\partial x}~\ln{[f(x)]}\right\}^2=4\,\int
dx\,\left[\psi'(x)\right]^2;\,\,\,\,\,(d\psi/dx= \psi').\ee Note
that for the uniform distribution $f(x)=constant$ one has $I=0$.
Focus attention now a system that is specified by a set of $M$
physical parameters $\mu_k$. We can write
$\mu_k = \langle A_{k}\rangle$ with $A_{k}= A_{k}(x).$
The set of $\mu_{k}$-values constitutes the prior knowledge. It
represents available empirical information. Let the pertinent
probability distribution function (PDF) be $f(x)$. Then, \be
\label{eq.1-4} \langle A_{k}\rangle\,=\,\int ~dx ~A_{k}(x) ~f(x),
\hspace{0.5cm} k=1,\dots ,M. \ee
  In this context it can be  shown (see for example
  \cite{pla7,reginatto}) that the {\it physically relevant}
  PDF $f(x)$  minimizes the FIM (\ref{eq.1-2}) subject to
  the prior conditions and the normalization condition.
In the celebrated MaxEnt approach of Jaynes' \cite{katz} one {\it
maximizes} the entropy, that behaves information-wise in opposite
fashion to that of Fisher's measure \cite{frieden3}. Normalization
entails $\int dx  f(x) = 1,$
 and, consequently, our Fisher-based extremization
problem adopts  the appearance \be \label{eq.1-6}\delta \left( I -
\alpha \int ~dx ~f(x) - \sum_{k=1}^M~\lambda_k\int ~dx ~A_{k}(x)
~f(x)\right) = ~0 \ee where we have introduced the $(M+1)$
Lagrange multipliers $\lambda_k$ ($\lambda_0=\alpha$).
In Ref. \cite{pla7} on can find
the details that lead from (\ref{eq.1-6}) to a Schr\"odinger's
equation (SE) that yields the desired PDF in terms of an amplitude
$\psi(x)$ referred to above before Eq. (\ref{eq.1-2}). This SE is
of the form

\be \label{eq.1-8} -~\frac{1}{2}~\frac{\partial^2 ~}{\partial x^2} \psi~-~\sum_{k=1}^{M}~
\frac{\lambda_{k}}{8}~ A_{k}\,\psi ~= ~ \frac{\alpha}{8}~ \psi, \ee
 which can be formally interpreted as the (real) Schr\"odinger equation
 for a particle of unit mass ($\hbar=1$) moving in the effective,
 ``information-related pseudo-potential"  \cite{pla7}

\be \label{eq.1-9} U~=~U(x)
=~-\frac{1}{8}~\sum_{k=1}^{M}\,\lambda_{k}~ A_{k}(x), \ee in which
the normalization-Lagrange multiplier ($\alpha /8$) plays the role
of an energy eigenvalue. The  $\lambda_k$ are fixed, of course, by
recourse to the available prior information. Note that $\psi(x)$
is always real in the case of one-dimensional scenarios, or for
the ground state of a real potential in N dimensions
\cite{richard}. In terms of the amplitudes $\psi(x)$ we have

\ben \label{eq.1-12-}  I & = & \,\int dx ~f \left(\frac{\partial
\ln{f} }{\partial x}\right)^2\,= \, \int dx ~ \psi_n^2 ~
\left(\frac{\partial \ln{\psi_n^2} }{\partial x} \right)^2\,=\, 4
\int d x ~ \left(\frac{\partial \psi_n }{\partial x} \right)^2
=\n \\
& =&\, -~4 \int \psi_n \frac{\partial^2 ~}{\partial x^2} \psi_n~dx
=\, -~4 \left\langle \frac{\partial^2 ~}{\partial x^2}
\right\rangle%
= \int ~ \psi_n \left(\alpha + \sum_{k=1}^M~\lambda_k~A_k\right)
\psi_n~dx,\n \een
\nd i.e.,\ben \label{eq.1-12} I=\,\alpha
 + \sum_{k=1}^M~\lambda_k\left\langle A_k\right\rangle. \een
 a  form that we will employ in our developments below.
 The connection between our variational solutions $f$  and thermodynamics
  was established in Refs. \cite{pla7} and \cite{flego} in the guise of
  reciprocity relations that  express  the
Legendre-transform structure of thermodynamics. They constitute
its essential formal ingredient \cite{deslog} and  were re-derived
\`a la Fisher in \cite{pla7}  by recasting  (\ref{eq.1-12}) in a
fashion that emphasizes the relevant independent variables

\ben \label{eq.1-13a} I(\left\langle
A_1\right\rangle,\ldots,\left\langle A_M\right\rangle) \,=\,\alpha
 + \sum_{k=1}^M~\lambda_k\left\langle
A_k\right\rangle. \een The Legendre transform changes the relevant
variables. In the case of $I$ this is

\be \label{eq.1-13b}\alpha= I(\left\langle
A_1\right\rangle,\ldots,\left\langle A_M\right\rangle) -
\sum_{k=1}^M~\lambda_k\left\langle A_k\right\rangle =
\alpha(\lambda_1,\ldots,\lambda_M), \ee so that we encounter the
three reciprocity relations  proved in \cite{pla7}

\be \label{RR-1} \frac{\partial \alpha}{\partial \lambda_{i}}= -
\langle A_i\rangle ~; \hspace{1.cm} \lambda_k \,=\, \frac{\partial
I }{\partial \left\langle A_k \right\rangle} ~ ;\hspace{1.cm}
\frac{\partial I}{\partial \lambda_{i}}=\sum_{k}^{M} \lambda_{k}
 \frac{\partial \langle A_{k}\rangle}{\partial \lambda_{i}},\ee
the last one being a generalized Fisher-Euler theorem.

\section{III. Present core-results}

\nd Our previous Fisher considerations lead to a scenario for
which

 \ben \label{SWE-1}
 H=-\frac{1}{2}~\frac{\partial^2~}{\partial x^2} ~+ U(x)~,
\hspace{0.5cm} U(x)= -~ \frac{1}{8} \sum_k \,\lambda_k\, A_k~,
\hspace{0.5cm} E_n=\frac{\alpha}{8}~.\een

\nd Enters here, as essential new ingredient in the present
considerations, the celebrated virial theorem \cite{virial} that
of course applies in this Schr\"odinger-scenario \cite{greiner}.
This theorem is intimately related to the reciprocity relations of
the preceding Section, as discussed in \cite{nuestro2} and states
that
 \ben \label{virial-4} \left\langle - ~ \frac{\partial^2
~}{\partial x^2}\right\rangle = \left\langle {x} ~ \frac{\partial
~}{\partial x} U({x})\right\rangle. \een  The potential function
$U(x)$ belongs to $\mathcal{L}_2$
 and thus admit of a series expansion in
$x,\,x^2,\,x^3,\,$etc. \cite{greiner}. The $A_k(x)$ themselves
belong to $\mathcal{L}_2$ as well and can be series-expanded in
similar fashion. This enables us to base our future considerations
on the assumption that the a priori knowledge refers to moments
$x^k$ of the independent variable, i.e.,
\be \langle A_k \rangle~=~ \langle x^k \rangle ~,  \ee and that
one possesses information on  $M$ moment-mean values
 $\langle x^k \rangle$. Our ``information" potential $U$   then reads
\be \label{virial-5} U(x)=  -~ \frac{1}{8} \sum_k \,\lambda_k\,
x^k\,.\ee and Eq. (\ref{virial-4})  allows one to immediately
obtain

\ben
\label{virial-6}
 \left\langle  \frac{\partial^2 ~}{\partial x^2}\right\rangle \,  =\, ~
\frac{1}{8}~\sum_{k=1}^{M}\, k \,\lambda_{k}~\left\langle
A_{k}\right\rangle;\hspace{1.2cm}(A_k=x^k),   \een and thus, via
 (\ref{virial-6}) and the above mentioned relation
 $I=-~4 \left\langle \partial_{xx}\right\rangle$,
 a useful, virial-related expression
for Fisher's information measure can be arrived at.

\ben \label{virial-7} I\, =\, -~ ~\sum_{k=1}^{M}\, \frac{k}{2}
\,\lambda_{k}~\langle x^{k}\rangle, \een which is an explicit
function of the M physical parameters $\langle x^{k}\rangle$ and
their respective Lagrange multipliers $\lambda_{k}$. Eq.
(\ref{virial-7}) encodes the information provided by the Virial
theorem. Thus, we have  two different ways of expressing $I$, namely,
(\ref{eq.1-12}) and (\ref{virial-7}). Interesting things happen if we
put them together.
Since $\lambda_k$ is given by (\ref{RR-1}) as $[\partial
I/\partial \langle x^k \rangle],$ inserting the reciprocity
relations (\ref{RR-1}) into (\ref{virial-7})  we are led to

\ben \label{gov-1} \frac{\partial I }{\partial \langle x^k
\rangle}\,=\,\lambda_k \hspace{0.5cm} &
\hspace{1.cm}\longrightarrow \hspace{1.cm} & I \,  =\,
-~\sum_{k=1}^{M}\,\frac{k}{2} ~\langle x^{k}\rangle ~
\frac{\partial I }{\partial \left\langle x^k \right\rangle}\,.
\een  Eq. (\ref{gov-1}) constitutes  an important result, since we
have now at our disposal  a {\it differential FIM-equation}.
Dealing with it should allow us to find $I$ in terms of the
$\langle x^k \rangle$ {\it without passing first through a
Schr\"odinger equation first}, a commendable achievement. This is
a linear partial differential equation that an extremal $I$ must
necessarily comply with. This constitutes one of the main present
results. It  is not clear (yet) whether from such an $I-$form we
can extract an amplitude $\psi$ satisfying a Schr\"odinger
equation. Our $I$ could however be related to an approximate
solution to Schr\"odinger's equation.
%%%%%%%%%%%%%%%%%%%%%%%%%%%%%%%%%%%%%%%%%%%%%%%%%%%%%%%%%%%%%%%%%%%%%%%%%%%%%%%%%
%%However, such issue is of no importance here. Why?
%%because our purpose is to codify the information provided by a set
%%of $M$ expectation values in an $I-$form. Previously, the only way
%%of doing that was to solve first a Schr\"odinger equation and {\it
%%now} we can bypass it.
For convenience we now recast our key
relations using dimensionless magnitudes
    \ben \label{gov-2}
    \forall~  \langle {A}_{k}\rangle \equiv \langle {x}^{k}\rangle ~\neq ~0 \hspace{0.2cm},\hspace{0.7cm}\mathcal{I}~=~    \frac{I}{[I]}~=~\frac{I}{[x]^2}~\hspace{0.2cm},\hspace{1.2cm}\langle\mathcal{X}_{k}\rangle ~=~ \frac{\langle x^k\rangle}{[\langle x^k\rangle ]} ~=~ \frac{\langle x^k\rangle}{[x]^{k}}~,\hspace{0.5cm}\een
 where $[I]$ and $[\langle x^k \rangle]$ denote the dimension of
$I$  and $\langle x^k \rangle$, respectively. Thus, the
differential equation that governs the FIM-behavior, i.e.,
(\ref{gov-1}),  can be translated into
 \ben  \label{gov-3}
 \mathcal{I}\,=\,-~\sum_{k=1}^{M}\,\frac{k}{2} ~\langle\mathcal{X}_{k}\rangle ~
\frac{\partial \mathcal{I} }{\partial\left\langle \mathcal{X}_k
\right\rangle} \,, \hspace{1.5cm}\mathcal{I}=
\mathcal{I}(\langle\mathcal{X}_{1}\rangle,\cdots,\langle\mathcal{X}_{M}\rangle),
\een \nd which is a first order linear nonhomogeneous equation
with $M$ independent variables.
All first order, linear partial differential equations (PDEs)  possess a solution that depends
 on an arbitrary function, called {\it the general solution} of the PDE.
In many physical situations this solution if {\it less important}
than other solutions called {\it complete ones}
\cite{pde,pde-1,pde-3}. Such complete solutions are particular PDE
solutions containing as many arbitrary constants as intervening
independent variables. Let us look now  for a special complete solution
of   our PDE (\ref{gov-3}), whose usefulness will be illustrated via two physical examples below
 (the treatment of general solutions is postponed to Section VI). We first set

  \ben \label{gov-4} \mathcal{I}=~\sum_{k=1}^M~\mathcal{I}_k ~ =
  ~\sum_{k=1}^M~ \exp{\left[~ g(\langle \mathcal{X}_k \rangle)~\right]},\een
   and substituting (\ref{gov-4}) into (\ref{gov-1}) leads to

\ben \label{gov-5} \mathcal{I} \,  =\,
-~\sum_{k=1}^{M}\,\frac{k}{2} ~\langle \mathcal{X}_{k}\rangle ~
g'(\langle \mathcal{X}_k \rangle)~\mathcal{I}_k\,.\een  The above
relation entails
 \ben \label{gov-6} ~
g'(\langle \mathcal{X}_k \rangle) = -~\frac{2}{k~\langle
\mathcal{X}_{k}\rangle} \hspace{1.cm}\longrightarrow \hspace{1.cm}
~ g(\langle \mathcal{X}_k \rangle) = -~\frac{2}{k}~\ln{\left|
\langle \mathcal{X}_{k} \rangle \right|}+c_k \, , \een where $c_k$
is an integration constant.  Finally, substituting (\ref{gov-6})
into (\ref{gov-4}) we arrive at

 \ben
\label{gov-7}
\mathcal{I}=\sum_{k=1}^{M}~C_k~\exp{\left(-~\frac{2}{k}~
\ln{\left|\langle \mathcal{X}_{k}\rangle \right|} \right)}~ ,
\hspace{1.2cm}C_k =~ e^{c_k}~>~0~,\een  which can be recast as
\ben \label{gov-8} \mathcal{I}(\langle \mathcal{X}_1 \rangle, ...
, \langle \mathcal{X}_M \rangle ) =~ \sum_{k=1}^{M}~C_k~ ~{ \left|
\langle \mathcal{X}_{k}\rangle \right|^{- {2}/{k}}}~ , \een or, in
function of the original input-quantities (\ref{gov-2})

\ben \label{gov-9} {I}(\langle {x}^1 \rangle, ... , \langle {x}^M
\rangle ) = ~ \sum_{k=1}^{M}~C_k~ ~{ \left| \langle {x}^{k}\rangle
\right|^{- {2}/{k}}}~,\een an intriguing  result.  We enumerate
below the main properties of this minimal $I$.

\begin{itemize}
    \item{\bf FIM-domain}

\nd    Obviously, it is
     \[{\it Dom}[I]=\left\{(\langle {x}^1 \rangle, ... , \langle {x}^M
\rangle) / \langle {x}^k \rangle~\in ~\Re_o \right\}\]

    \item{\bf FIM-monotonicity}

Differentiating (\ref{gov-9}) we obtain
\ben \label{prop-2}  \frac{\partial I}{\partial \langle
x^k \rangle} =~-~ \frac{2}{k~\langle x^k \rangle}I_k~=~-~ \frac{2}{k \, \langle x^k \rangle} \, C_k
~\left| \langle {x}^{k}\rangle \right|^{- {2}/{k}}~~~ . \een
Therefore, if $
\langle x^k \rangle ~ > ~ 0 ~$ , $I$ is a monotonically decreasing
function in the $ \langle x^k \rangle$-direction.  Also, for $ \langle x^k \rangle ~ >
~ 0 ~$, from the reciprocity relations (\ref{RR-1}) we have, \ben
\label{prop-2-RR}
 \lambda_k ~=~-~ \frac{2}{k}~C_k \, \langle x^k \rangle^{-(2+k)/k} \,
<~0. \een

    \item{\bf FIM-convexity}

\nd This is a necessary property, since the entropy is concave. By
differentiation of the expression (\ref{prop-2}) one obtains \ben
 \label{prop-3} \frac{\partial^2 I }{\partial
\langle x^n \rangle\partial \langle x^k \rangle}
   \,=\, \left(1+ \frac{k}{2} \right)  \frac{4}{k^2}  ~
 ~ C_k~\left| \langle {x}^{k}\rangle \right|^{- {2}(1+k)/{k}}~\delta_{kn},\een
 from which we can assert that the Fisher measure is  a convex
function. It is then guaranteed that the inverse of $\partial_k
\partial_j \bar{\alpha}$ exists.

\end{itemize}

\section{IV. The reference quantities $C_k$}

\nd FIM is an estimation measure known to obey the Cramer Rao-bound (\ref{rao})
 \cite{frieden2}. The best estimator exhibits a CR relation as close
to unity as possible. Thus, the reference quantities $C_k$ should be chosen
in a manner that respects this condition. Here we  are interested in simple
situations that illustrate the concomitant procedure. More involved situations will be treated elsewhere.

\nd Since the reference quantities
$C_k$  contain important information concerning the reference
system with respect the which prior conditions are experimentally determined,
it is convenient to start by choosing an appropriate reference one.

\subsection{Minimum of the information potential}

\nd We consider it reasonable to incorporate at the outset, within the
$I-$form, information concerning the minimum of the information
potential $U(x)$. Assume that this information-potential  \ben
U(x)=-~\frac{1}{8}\sum_{x=1}^M \lambda_k x^k, \n \een achieves its
absolute minimum at the ``critical point" $x=\xi$, \ben
U^{~'}(\xi)~=~0~, \hspace{2.cm} U_{min}~=~U(\xi).\een \nd Effecting the FIM-translational transform $u=x-\xi $ leads us to
\ben \label{sys-1} I=-~ \sum_{k=1}^{M}~\frac{k}{2}~\lambda_k~
\langle {x}^{k}\rangle = ~-~
\sum_{k=1}^{M}~\frac{k}{2}~\lambda_k^*~ \langle {u}^{k}\rangle',
\een with (see the Appendix) \ben
 \label{sys-2} \lambda_k^*~= ~-~\frac{8}{k!}~U^{(k)}(\xi)~, \hspace{1.2cm}
  \langle {u}^{k}\rangle' \,= \,\langle (x-\xi)^{k}\rangle~  \een
where $U^{(k)}(\xi)$ is the $k^{th}$ derivative of U(x) evaluated
at $x=\xi$ and $\langle ~ \rangle'$ indicates that the pertinent mean value ($x-$moment) is
evaluated for  translation-transformed eigenfunctions. \nd The
corresponding FIM-explicit functional expression is built up with
the $N-$non-vanishing momenta ($N < M$) ($\langle u^k \rangle'
\neq 0$) and is given by
\ben \label{sys-3} {I} = ~
\sum_{k=2}^{N}~C_k ~{ \left| \langle {u}^{k}\rangle' \right|^{-2/k}}~= ~ \sum_{k=2}^{N}~C_k ~{ \left| \langle (x-\xi)^{k}\rangle
\right|^{- 2/k}}~,\een \nd where we kept in mind that
$\lambda_1^*= - 8 U'(\xi)~=~0.$  A glance at the above
expression suggests that we  re-arrange things in the fashion
\ben \label{CR-1} {I} = ~C_2~\left| \langle (x-\xi)^{2}\rangle
\right|^{- 1}+ \sum_{k=3}^{N}~C_k ~{ \left| \langle
(x-\xi)^{k}\rangle \right|^{- 2/k}}~.\een Taking now into account
that

\ben \label{CR-2} \left\{ \begin{array}{l}
\langle~ x-\xi~\rangle =0\\
\langle (x-\xi)^{2}\rangle =\langle x^{2}\rangle-2 \xi \langle x\rangle+\xi^2\\
\end{array} \right. \hspace{0.7cm}\longrightarrow\hspace{0.7cm}
 \left\{\begin{array}{l}
\langle~ x~\rangle ~=~\xi\\
\langle (x-\xi)^{2}\rangle =\langle x^{2}\rangle- \langle
x\rangle^2=\sigma^2
\end{array} \right. \hspace{0.7cm}
~\een  we get \ben \label{CR-3} {I} =   ~\frac{C_2}{\sigma^2}~+
\sum_{k=3}^{N}~C_k ~{ \left| \langle (x-\xi)^{k}\rangle \right|^{-
2/k}}~,\een from which we obtain \ben \label{CR-4} {I}~\sigma^2 =
~C_2~+ ~\sigma^2 \sum_{k=3}^{N}~C_k ~{ \left| \langle
(x-\xi)^{k}\rangle \right|^{- 2/k}}~\geq ~1~.\een Therefore, $I$
preserves the well-known Cramer-Rao $I-$bound \cite{frieden3} $
I~{\sigma^2}~\geq~ 1.$ The above  seems to indicate that if no moment $k \ge3$  is a priori known, the lower
bound can be reached for $C_2=1$.
For $k\ge 3$ additional considerations apply that will be
discussed elsewhere.

\section{V. Two physical examples}

\nd So as to illustrate  the above considerations we are going to
consider two simple and instructive examples. We take the mass
$m=1$ and $\hbar=1$.

\subsection{Harmonic oscillator (HO)}

 \nd The  prior information is given by
 \ben \label{ex-2PK}  \langle x^2\rangle ~ = ~\frac{1}{2\omega}~, \hspace{1.5cm} M~=~1~, \hspace{1.5cm}k~=~2~.\een

\nd The minimum of the potential function obtains at the origin  $\xi = 0$,
\[U(x)=-\frac{1}{8}\lambda_2~x^2 \hspace{1.cm}\longrightarrow\hspace{1.cm}U'(\xi)=-\frac{1}{4}\lambda_2~\xi=0
\hspace{1.cm}\longrightarrow\hspace{1.cm}\xi=0.\]

\nd The pertinent FIM  can be obtained using (\ref{sys-3}) with $u=x-\xi=x$,
\ben  I ~ = ~I({\langle x^2\rangle})~= ~C_2~
\langle x^2 \rangle ^{- 1} ~,\n \een
and, the CR bound is saturated when $C_2=1$,
\ben \label{ex-2a} I~\langle x^2 \rangle ~=~C_2~=~1 \hspace{1.cm}\Longrightarrow \hspace{1.cm}I ~ =  ~\langle x^2 \rangle ^{- 1} ~.  \een

\nd The corresponding Lagrange multiplier can be obtained by recourse to the reciprocity relations (\ref{RR-1}) and (\ref{ex-2a}),
\ben \label{ex-2b} \lambda_2 \,=~ \frac{\partial I }{\partial \langle x^2 \rangle}
\,=\, -~ \langle x^2 \rangle ^{- 2}~. \een

\nd The prior-knowledge (\ref{ex-2PK}) is encoded into the FIM (\ref{ex-2a}),
and the Lagrange multiplier $\lambda_2$ (\ref{ex-2b}),
\ben \label{ex-2c} I ~ = ~\langle x^2\rangle~^{-1}= ~ 2 \omega \,;
 \hspace{1.5cm}  \lambda_2 ~=~-~ \langle x^2 \rangle
~^{- 2}~= ~ - ~ 4\omega^2~. \een
and the $\alpha -$value can be obtained from (\ref{eq.1-13b}),
\ben \label{ex-2d} \alpha  ~=~ I-~\lambda_2 ~ \langle x^2\rangle
=~4~\omega, \een
as we expect.

\subsection{Harmonic oscillator in a uniform field}

\nd We consider a charged unit-mass particle moving in the HO
 potential. The electrical charge is q and there is a uniform electric field $\epsilon$,  in the $x-$direction.
Our prior knowledge is given by \cite{greiner}
\ben \label{ex-3}   \langle x \rangle ~ =  ~
\frac{q~\epsilon}{\omega^2}~, \hspace{2.cm} \langle x^2\rangle ~ =
~ \frac{1}{2\omega} + \left(\frac{q~\epsilon}{\omega^2}\right)^2~.\een

\vskip 4mm

\nd We look first for the $\xi$-point at which $U(x)$ is minimal.
    \[U(x)=-\frac{1}{8}\left(\lambda_1~x+\lambda_2~x^2 \right)\]
\ben \label{ex3-F1}
U'(\xi)~=~-\frac{1}{8}\left(\lambda_1+2\lambda_2~\xi \right)~=~0 \hspace{1.cm}\longrightarrow \hspace{1.cm}  \xi=-~\frac{\lambda_1}{2~\lambda_2}\hspace{0.3cm}~. \een
\nd The translational transform $u=x-\xi $ implies that
\ben \label{ex3-F2} \langle u\rangle'~=~\langle x-\xi \rangle =
\langle x \rangle - \xi~,\hspace{1.2cm}
 \langle u^2\rangle'~=~\langle (x-\xi)^{2}\rangle =
\langle x^{2}\rangle - 2 \xi \langle x \rangle + \xi^2~.\een

\nd The translation transformed FIM is now given by
\ben \label{ex3-F3} I~=~ C_2~ \langle
{u}^{2}\rangle^{'-1}.\een

\nd and, the CR bound is saturated when $C_2=1$,
\ben \label{ex3-F4} I~\langle u^2 \rangle' ~=~C_2~=~1 \hspace{1.cm}\Longrightarrow \hspace{1.cm} I ~ =  ~ \langle u^2 \rangle^{'- 1} ~.  \een
 The reciprocity relations lead us to
\ben \label{ex3-F5}
\lambda_1 \,&=&\,\frac{\partial~I}{\partial \langle x \rangle}\,=\,
\,\frac{\partial~I}{\partial \langle u^2 \rangle'}\,
\frac{\partial \langle u^2 \rangle'}{\partial\langle x\rangle}~= ~
-\langle {u}^{2}\rangle^{'-2}~(-2~\xi)~ \\
\label{ex3-F6}
 \lambda_2 \,&=&\,\frac{\partial~I}{\partial \langle x^2 \rangle}\,=\,
\,\frac{\partial~I}{\partial \langle u^2 \rangle'}\,
\frac{\partial \langle u^2 \rangle'}{\partial \langle
x^2\rangle}~= ~ -~ \langle {u}^{2}\rangle^{'-2}~. \een
From the prior knowledge (\ref{ex-3}) and using (\ref{ex3-F2}) we have
\ben \label{ex3-F7} \langle x \rangle =
\xi~=~\frac{q~\epsilon}{\omega^2}~,\een \ben \label{ex3-F8}
\langle u^2\rangle'~=~ \langle x^{2}\rangle -
\xi^2~=~\frac{1}{2\omega} +
\left(\frac{q~\epsilon}{\omega^2}\right)^2~
-~\left(\frac{q~\epsilon}{\omega^2}\right)^2~=~\frac{1}{2\omega}~,\een
then, inserting (\ref{ex3-F7}) and (\ref{ex3-F8}) into (\ref{ex3-F4}) -
(\ref{ex3-F6}) we get
\ben \label{ex3-F9} I~=~ \langle
{u}^{2}\rangle^{'-1}~=~
\left(\frac{1}{2\omega}\right)^{-1}~=~2\omega,\een  \ben
\label{ex3-F10}
 \lambda_1 \,&=&\,2~\xi~\langle {u}^{2}\rangle^{'-2}~=
  ~2~\frac{q~\epsilon}{\omega^2}~(2~ \omega)^2~ =~8~q~\epsilon \\
\label{ex3-F11}
 \lambda_2 \,&=&\,-\langle {u}^{2}\rangle^{'-2}~=~- (2~ \omega)^2~=~- 4~ \omega^2
\een The corresponding translational transform $\bar{\alpha} -$value can be obtained substituting (\ref{ex3-F9})-(\ref{ex3-F11}) into (\ref{eq.1-13b}),
 \ben \label{ex3-F12} \bar{\alpha}  = \, I~-\lambda_1 \langle
x\rangle-\lambda_2 \langle x^2\rangle~=~2\omega -8 q \epsilon
\frac{q \epsilon}{\omega^2} ~+ 4 \omega^2 \left(\frac{1}{2\omega}
+ \left(\frac{q~\epsilon}{\omega^2}\right)^2\right)~ =~4\omega~,
\een and the corresponding ${\alpha} -$value is given by (see
Appendix),
 \ben \label{ex3-F13} \alpha=\bar{\alpha}+8~U(\xi) ~=~4\omega~-4~\frac{q^2\epsilon^2}{\omega^2},
\een
as we expect.

\section{VI. General Solution of the differential FIM-equation}

\nd We discuss here this issue for the sake of completeness. Our
$FIM-$equation is a first order linear nonhomogeneous differential
equation. We are following \cite{pde,pde-1,pde-3} in looking for
the general solution. For a first-order PDE, the {\it method of
characteristics}
 allows one to encounter useful curves (called characteristic curves or just
characteristics) along which the PDE becomes an ordinary
differential equation (ODE). Once the ODE is found, it can be
solved along the characteristic curves and transformed into a
solution for the original PDE.

\nd The  characteristic system of Eq. (\ref{gov-3}) is  \ben
{\label{pde-I2}}
-~\frac{d\langle\mathcal{X}_{i}\rangle}{(i/2)\langle\mathcal{X}_{i}\rangle
}=-~\frac{d\langle\mathcal{X}_{j}\rangle}{(j/2)\langle\mathcal{X}_{j}\rangle}
=\frac{d\mathcal{I}}{\mathcal{I}}~,\hspace{1.cm}
i,j=1,\cdots,M,\een leads (for $\langle\mathcal{X}_{1}\rangle\neq
0$) to

\ben {\label{pde-I3}} \frac{d\langle\mathcal{X}_{k}\rangle}{(k/2)\langle\mathcal{X}_{k}\rangle}
=\frac{d\langle\mathcal{X}_{1}\rangle}{(1/2)\langle\mathcal{X}_{1}\rangle
} \hspace{0.5cm}\longrightarrow\hspace{0.9cm}
\frac{2}{k}~\ln{\left|\langle\mathcal{X}_{k}\rangle\right|}+{c_k}&=&
2~\ln{\left|\langle\mathcal{X}_{1}\rangle\right|}+{c_1}\n\\
\ln{~\left[e^{c_k}~\left|\langle\mathcal{X}_{k}\rangle\right|^{2/k}\right]}&=&\ln{\left[e^{c_1}~\left|\langle\mathcal{X}_{1}\rangle\right|^{2}\right]}\n\\
~e^{c_k}~\left|\langle\mathcal{X}_{k}\rangle\right|^{2/k}&=&e^{c_1}~\left|\langle\mathcal{X}_{1}\rangle\right|^{2}\n\\
&\downarrow & \n \\
b_{k-1}\equiv
e^{c_k-c_1}&=&~{\left|\langle\mathcal{X}_{1}\rangle\right|^{2}}
{\left|\langle\mathcal{X}_{k}\rangle\right|^{-2/k}~}\\
\n\\{\label{pde-I4}}
\frac{d\mathcal{I}}{\mathcal{I}}= -\frac{d\langle\mathcal{X}_{1}\rangle}{(1/2)
\langle\mathcal{X}_{1}\rangle}
\hspace{1.5cm}\longrightarrow\hspace{0.9cm}
\ln{\left|\mathcal{I}\right|}+{c_{\mathcal{I}}}~&=&~
- 2~\ln{\left|\langle\mathcal{X}_{1}\rangle\right|}+{c_1}\n\\
\ln{~\left[e^{c_{\mathcal{I}}}~\left|{\mathcal{I}}\right|\right]}&=& \ln{\left[e^{c_1}~\left|\langle\mathcal{X}_{1}\rangle\right|^{-2}\right]}\n\\
e^{c_{\mathcal{I}}}~\left|{\mathcal{I}}\right|&=&
e^{c_1}~\left|\langle\mathcal{X}_{1}\rangle\right|^{-2}\n\\
&\downarrow & \n \\
b_M\equiv e^{c_1-c_{\mathcal{I}}}&=&~{\left|\langle\mathcal{X}_{1}\rangle\right|^{2}~}
{\left|\mathcal{I}\right|}~.\een
\nd We have now constructed an integral basis for the
characteristic system (\ref{pde-I2}) \ben {\label{pde-I5}}
b_1=u_1(\langle\mathcal{X}_{1}\rangle,...,\langle\mathcal{X}_{M}\rangle,\mathcal{I})~,  ~.~.~.~,~
b_M=u_{M}(\langle\mathcal{X}_{1}\rangle,...,\langle\mathcal{X}_{M}\rangle,\mathcal{I})~, \een \nd and the
general solution of equation (\ref{gov-3}) defined as
\ben {\label{pde-s6}} \Phi(u_1,u_2,~.~.~.~,u_{M})~=~0, \een is
given by
\ben {\label{pde-s7}}
\Phi\left({\left|\langle\mathcal{X}_{1}\rangle\right|^{2}}
{\left|\langle\mathcal{X}_{2}\rangle\right|^{-1}},
\cdots,{\left|\langle\mathcal{X}_{1}\rangle\right|^{2}}{\left|\langle\mathcal{X}_{k}\rangle\right|^{-2/k}~},\cdots,{\left|\langle\mathcal{X}_{1}\rangle\right|^{2}}{\left|\langle\mathcal{X}_{M}\rangle\right|^{-2/M}~},
{\left|\langle\mathcal{X}_{1}\rangle\right|^{2}~}
{\left|\mathcal{I} \right|}\right)~=~0,\n
\een \nd where $\Phi$ is an arbitrary function of the $M$
variables. Solving this equation for $\mathcal{I}$ yields a
solution of the  explicit form
\ben {\label{pde-I8}} \mathcal{I}
~=~{~~\left|\langle\mathcal{X}_{1}\rangle\right|^{-2}~}~
\Psi \left(
{\left|\langle\mathcal{X}_{1}\rangle\right|^{2}}
{\left|\langle\mathcal{X}_{2}\rangle\right|^{-1}},
\cdots,{\left|\langle\mathcal{X}_{1}\rangle\right|^{2}}{\left|\langle\mathcal{X}_{k}\rangle\right|^{-2/k}~},\cdots,{\left|\langle\mathcal{X}_{1}\rangle\right|^{2}}{\left|\langle\mathcal{X}_{M}\rangle\right|^{-2/M}~}
\right),
\een \nd where $\Psi$ is an arbitrary function of ($M-1$)
variables.

\vspace{0.5cm}

\nd {\bf Cauchy problem and the existence and uniqueness of the solution to our PDE}

\vspace{0.3cm}

\nd One of the fundamental aspects so as to have a useful PDE for modeling physical  systems  revolves around the existence and  uniqueness of the solutions to
the Cauchy problem. Here we show that such  requirements are satisfied by our pertinent solutions.
\nd We start by casting (\ref{gov-3}) in the normal form

\ben {\label{pde-s14}}
 \frac{\partial{\mathcal{I}}}{\partial \langle\mathcal{X}_{1}\rangle} =F\left(\langle\mathcal{X}_{1}\rangle,\cdots,\langle\mathcal{X}_{M}\rangle,
 {\mathcal{I}}, \frac{\partial{\mathcal{I}}}{\partial \langle\mathcal{X}_{2}\rangle},\cdots, \frac{\partial{\mathcal{I}}}{\partial \langle\mathcal{X}_{M}\rangle}\right)
\een
 where \ben
{\label{pde-I15}}F_{\mathcal{I}}=F\left(\langle\mathcal{X}_{1}\rangle,
\cdots,\langle\mathcal{X}_{M}\rangle,
 {\mathcal{I}}, \frac{\partial{\mathcal{I}}}{\partial \langle\mathcal{X}_{2}\rangle},\cdots, \frac{\partial{\mathcal{I}}}{\partial \langle\mathcal{X}_{M}\rangle}\right)=
 -\frac{2}{\langle\mathcal{X}_{1}\rangle}\left[\mathcal{I}+
~\sum_{k=2}^{M}\,\frac{k}{2} ~\langle\mathcal{X}_{k}\rangle
\frac{\partial \mathcal{I} }{\partial\left\langle \mathcal{X}_k
\right\rangle}
\right],~ \een
\nd and we see that $F$ is a real function of class $C^2$ in a
neighborhood of \ben {\label{pde-s16}} \langle\mathcal{X}_{1}\rangle=a~,
\hspace{0.2cm}\langle\mathcal{X}_{k}\rangle=\xi_{k-1}~, \hspace{0.2cm}
\mathcal{I}({\xi_1,...,\xi_{M-1}})=c~, \hspace{0.2cm}
\left.\frac{\partial{\mathcal{I}}}{\partial
\langle\mathcal{X}_{k}\rangle} \right|_{\xi_1,...,\xi_{M-1}}=d_{k-1}~,
\hspace{0.3cm}k=2,\cdots,M \hspace{0.3cm} \n \een \nd Then, if
$\psi(\langle\mathcal{X}_{2}\rangle,...,\langle\mathcal{X}_{M}\rangle)~$ is also a function of class $C^2$ such that
 \ben \label{pde-I17}
\psi(\xi_1,...,\xi_{M-1})=c~,\hspace{1.5cm}
\left.\frac{\partial{\psi}}{\partial
\langle\mathcal{X}_{k}\rangle}\right|_{\xi_1,...,\xi_{M-1}}=d_{k-1}~,
\hspace{0,5cm}k=2,...,M. \een \nd exists a solution
$\mathcal{I}$ of (\ref{pde-s14}) in a neighborhood of
$\langle\mathcal{X}_{1}\rangle=a$ and $\langle\mathcal{X}_{k}\rangle=\xi_{k-1}$
 that satisfies \ben
\label{pde-s18} \mathcal{I}(a,
\langle\mathcal{X}_{2}\rangle,\cdots,\langle\mathcal{X}_{M}\rangle)=\psi(\langle\mathcal{X}_{2}\rangle,\cdots,\langle\mathcal{X}_{M}\rangle) \een
\nd and is of class $C^2$.

\nd Regarding Cauchy-uniqueness,  it is known that if $F$
satisfies the Lipschitz condition \cite{pde-6}, \ben
\label{pde-I19}
\left|F_{\mathcal{I}^{~'}}-F_{\mathcal{I}}\right|~\leq ~K_1
\sum_{k=2}^M \left| \frac{\partial{\mathcal{I}^{~'}}}{\partial
\langle\mathcal{X}_{k}\rangle}-
\frac{\partial{\mathcal{I}}}{\partial
\langle\mathcal{X}_{k}\rangle}\right| +K_2
\left|\mathcal{I}^{~'}-\mathcal{I}\right|~ \hspace{1.cm}K_1,~K_2~=
const.\hspace*{1.cm} \een
\[F_{\mathcal{I}^{~'}}=F\left(\langle\mathcal{X}_{1}\rangle,
\cdots,\langle\mathcal{X}_{M}\rangle,
 {\mathcal{I}^{~'}}, \frac{\partial{\mathcal{I}^{~'}}}{\partial \langle\mathcal{X}_{2}\rangle},\cdots, \frac{\partial{\mathcal{I}^{~'}}}{\partial \langle\mathcal{X}_{M}\rangle}\right)\]
 then, the solution of the initial value
problem for (\ref{pde-s14}) is unique. Note that in our case the
above condition  is verified always since the Legendre structure
the theory guarantee that

\ben \label{pde-s20}
F_{\mathcal{I}}=
\frac{\partial{\mathcal{I}}}{\partial \langle\mathcal{X}_{1}\rangle}\propto
\frac{\partial{I}}{\partial \langle x\rangle}= \lambda_1 < \infty. \een

%%%%%%%%%%%%%%%%%%%%%%%%%%%%%%%%%%%%%%%%%%%%

\section{VII. Conclusions}

\nd It was known for some time  \cite{pla7} that, for
one-dimensional scenarios,  a minimal Fisher's information measure
$I$  is associated  to a Schr\"odinger's equation in which the
role of the potential is played by a weighted sum of the a priori
known expectation values of $M$ quantities $\langle x^{k}\rangle$.
The weights are the corresponding Lagrange multipliers
$\lambda_k$. Legendre reciprocity relations linking $I$, the
$\lambda_k$ and the $\langle x^{k}\rangle$ apply. \vskip 4mm \nd
We have here demonstrated that insertion of  virial theorem-tenets
into this Legendre structure  leads to a differential equation for
$I$. The equation is analytically solvable and its solution
provides us with  explicit new expressions for $I$ in terms of the
input-information  contained in the $M$  expectation values
$\langle x^{k}\rangle$. In other words, we can directly codify the
information provided by such set of expectation values in an
$I-$form without previous appeal to a Schr\"odinger equation.
 Additionally, this (partial) differential equation (PDE) may be viewed as a new constraint
 that the solutions $\psi_n$  of Schr\"odinger's equation (SE) must satisfy. In
 this regard, our PDE may lead to useful criteria for checking
 approximate treatments of SE.  If $\phi_n$ is an approximate
 SE-solution, the fact that the functional $I[\phi_n]$  verifies  the
 PDE would constitute an indicator of the quality of the approximate solution $\phi_n$.
\vskip 4mm \nd An application to
 simple examples  has illustrated these
considerations. Of course, as is the case  in the MaxEnt
environment, the usefulness of (\ref{gov-9}) depends on how
adequate is our input information for describing the situation at
hand. Maximal entropy or minimum FIM are just the best ways to
exploit  that knowledge. \vskip 4mm

%%\nd It is important to note that, once in possession
%% of this minimal $I$, one does not need to ever address any explicit
%% $I-$minimization task nor  solve the FIM-associated Schr\"odinger equation, a fact that should be of
 %%%interest to the large number of Fisher practitioners. Additionally,
 %%some new light on the Legendre transform substructure that
 %%%underlies Schr\"odinger's equation has been here shed.

\vspace{1.cm}
\appendix

\section{Appendix: Translation transform of FIM}

\nd The potential function \ben \label{ap-1}
U(x)=-~\frac{1}{8}\sum_{k=1}^M \lambda_k x^k. \n \een

\nd can be Taylor-expanded about $x=\xi$  \ben \label{ap-2}
U(x)=\sum_{k=0}^M \frac{U^{(k)}(\xi)}{k~!} (x-\xi)^k \n \een

\vspace{0.2cm}

\nd The translational transform $u=x-\xi $  leads to \ben
\label{ap-3} \bar{U}(u)=~{U}(u+\xi)=~
\sum_{k=0}^{M}~\frac{U^{(k)}(\xi)}{k~!} u^k,\een \nd which can be
recast as \ben \label{ap-4} \bar{U}(u)=-~\frac{1}{8}
\sum_{k=0}^{M}~\lambda^{*}_k u^k ~,\een \nd with \ben \label{ap-5}
\lambda^{*}_k ~\equiv ~ -~ 8 ~\frac{U^{(k)}(\xi)}{k~!}= -
~\frac{8}{k!}\sum_{j=1}^{M}~j(j-1)(j-2)\cdots(j-k+1)~\lambda_j~
{\xi}^{j-k}~. \een

\nd The FIM-translational transform $u=x-\xi $ is obtained from
(\ref{eq.1-12}) in the fashion ($\langle ~ \rangle'$ indicates
that the moment is calculated for translation-transform
eigenfunctions)

\ben \label{ap-6} I \, =\, -~4 \int {\psi} \frac{\partial^2
~}{\partial x^2} {\psi}~dx =\, -~4 \int \bar{\psi}
\frac{\partial^2 ~}{\partial u^2} \bar{\psi}~du =\, -~4
\left\langle \frac{\partial^2 ~}{\partial u^2} \right\rangle',
\een  where $\bar{\psi}=\bar{\psi}(u)$ is the translation
transform of ${\psi}(x)$. Now, using  the translation transform of
(\ref{eq.1-8}) one easily finds
 \ben \label{ap-7} I\,=\, \int ~ \bar{\psi}_n \left(\alpha +
\sum_{k=0}^M~\lambda^{*}_k~u^k\right) \bar{\psi}_n~du ~, \een \nd
and one realizes that

 \ben \label{ap-8} I \,=\,\alpha
 + \sum_{k=0}^M~\lambda^{*}_k  \langle u^k\rangle' \,=\,\bar{\alpha}
 + \sum_{k=1}^M~\lambda^{*}_k  \langle u^k\rangle'~, \een
where \ben \label{ap-9}
\bar{\alpha}=\alpha+\lambda^{*}_0=\alpha-8U(\xi). \een

\vspace{0.2cm}

\nd Also, the virial theorem (\ref{virial-4}) leads to

\ben \label{ap-10}I= 4~\left\langle ~ \frac{\partial^2
~}{\partial u^2}\right\rangle' = -~4~\left\langle {u} ~ \frac{\partial
~}{\partial u} \bar{U}(u)\right\rangle' ~
=~-~ \sum_{k=1}^{M}~\frac{k}{2}~\lambda^{*}_k~ \langle {u}^{k}\rangle' \een

\vspace{0.5cm}

\nd The translation-transformed moments $\langle {u}^{k}\rangle'$
are related to the original moments as

\ben \label{ap-11} \langle {u}^{k}\rangle' \, =\,  \int {u}^{k}~\bar{\psi}^2(u) ~du
\, =\,  \int {u}^{k}~{\psi}^2(u+\xi) ~du \,= \,  \int (x-\xi)^{k}~{\psi}^2(x) ~dx  \,= \,\langle (x-\xi)^{k}\rangle~\n  \een

\vspace{0.5cm}

\nd By recourse to the Newton-binomial  we write \ben
\label{ap-12} \int (x-\xi)^{k}~{\psi}^2(x) ~dx
~=~\sum_{j=1}^k{~(-1)^j~ \left( \begin{array}{c} k\\j  \end{array}
\right)~\xi^j~\int x^{k-j}~{\psi}^2(x) ~dx}, \een \nd and then we
finally have
    \ben \label{ap-13}
    \langle u^k \rangle' ~ = \,\langle (x-\xi)^{k}\rangle~ = ~\sum_{j=1}^k{~(-1)^j~
\left( \begin{array}{c} k\\j    \end{array} \right)~\xi^j~\langle
x^{k-j}\rangle}.\een

\vspace{0.5cm} \nd {\bf Acknowledgment:}  Partial support from the
programs FQM-2445 and FQM-207 of the Junta de Andalucia-Spain and
from CONICET (Argentine Agency) is acknowledged.

%%%%%%%%%%%%%%%%%%%%%%%%%%%%%%%%%%%%%%%%%%%%%%%%%%%%%%%%%%%%%%%%%%%%%%%%


\begin{thebibliography}{99}
\bibitem{virial} For any quantum system in stationary state,
 with a Hamiltonian does not involve time explicitly,
$H=- \frac{\hbar^2}{2m}~\vec{\nabla}^2~+~U(\vec{x}),~$the virial
theorem states  that \cite{greiner} $ \left\langle  -
\frac{\hbar^2}{m}\vec{\nabla}^2 \right\rangle = \left\langle
\vec{x}.\vec{\nabla} U(\vec{x})\right\rangle,$ where the
expectation value is taken for stationary states of the
Hamiltonian.
\end{thebibliography}
\end{document}